\documentclass[12pt, twoside]{amsart}
\usepackage[english]{babel}
\theoremstyle{plain}
\addtocounter{page}{1}
\usepackage{psfrag}
\usepackage{graphicx}
\newtheorem{Theorem}{Theorem}[section]

\newtheorem{Proposition}{Proposition}[section]

\newcommand{\be}{\begin{equation}}
\newcommand{\ee}{\end{equation}}
 
\newcommand{\rec}[1]{{(\ref{#1})}} 
\newcommand{\lapfr}{(-\Delta)^\frac{1}{2}} 
\newcommand{\ba}{\begin{array}} 
\newcommand{\ea}{\end{array}}

 \def\11{1\!\!1}

\def\p{\partial}
\def\ti{\tilde}

\def\ds{\displaystyle}

 \newcommand{\R}{\mathbb{R}}

\newcommand{\C}{\mathbb{C}}

\addtocounter{page}{1}

\def \R {{\mathbb {R}}}

\def \C {{\mathbb C}}

\title[Pohozaev-Type Identities]{Some  Remarks on  Pohozaev-Type Identities}
 \author{Francesca Da Lio }
\address{Department of Mathematics, ETH Z\"urich, R\"amistrasse 101, 8092 Z\"urich, Switzerland}
\email{fdalio@math.ethz.ch }

\begin{document}

\begin{abstract}
{The aim of this note is to discuss in  more detail the Pohozaev-type identities that have been
 recently established by the author, P. Laurain and T. Rivi\`ere  in  \cite{DLLR} in the framework of half-harmonic maps defined   either on $\R$ or on the sphere $S^1$ with values into a closed manifold ${\mathcal{N}}^n\subset \R^m$.  Weak  half-harmonic maps  are critical points of the following nonlocal energy
  \begin{equation}\label{fracenergyab}
{\mathcal{L}}_{\R}^{1/2}(u):=\int_{\R}|(-\Delta)^{1/4}u|^2\ dx~~\mbox{or}~~{\mathcal{L}}_{S^1}^{1/2}(u):=\int_{S^1}|(-\Delta)^{1/4}u|^2\ d\theta.
\end{equation}
If $u$ is a sufficiently  smooth   critical point of \eqref{fracenergyab}
then it satisfies the following   equation  of stationarity
\begin{equation}\label{statequ}
\frac{du}{dx}\cdot (-\Delta)^{1/2} u=0~~\mbox{a.e in $\R$}~~\mbox{or}~~\frac{\partial u}{\partial \theta}\cdot (-\Delta)^{1/2} u=0~~\mbox{a.e in $S^1$.}\end{equation}
As it was announced  in \cite{DLLR}, by using the  invariance of \eqref{statequ} in $S^1$   with respect to the trace of the M\"obius transformations of the $2$ dimensional disk we derive a countable family of relations  involving
the Fourier coefficients of weak half-harmonic maps $u\colon S^1\to {\mathcal{N}}^n.$ In the same spirit we also provide   as many  Pohozaev-type identities in $2$-D  for stationary
harmonic maps as conformal vector fields in $\R^2$ generated by holomorphic functions.}\par
\medskip \noindent {\sc{2010 MSC.}}  35R11,58E20, 30C20, 42A16.

 \noindent {\sc{Keywords.}} Fractional harmonic maps, harmonic maps, M\"obius transformations, Pohozaev formulas, Fourier coefficents.

 \end{abstract}

\maketitle
\section{
Introduction}
The notion of weak $1/2$-harmonic maps   into a $n$-dimensional closed (compact without boundary) manifold ${\mathcal{N}}^n\subset \R^m$ has been introduced by Tristan Rivi\`ere  and the author in
\cite{DLR1,DLR2}.  Since then the theory of fractional harmonic maps have received a lot of attention in view of their application to important geometrical problems (see e.g \cite{DL1,DLSR} for an overview of the theory).\par

These maps
are critical points of the fractional energy on ${\R}^k$
\begin{equation}\label{fracenergy}
{\mathcal{L}}^{1/2}(u):=\int_{{\R}^k}|(-\Delta)^{1/4}u|^2\ dx^k
\end{equation}
 within
\[
\dot H^{1/2}({\R}^k,{\mathcal{N}}^n):=\left\{ u\in \dot H^{1/2}({\R}^k,{\R}^m)\ ;\ u(x)\in {\mathcal{N}}^n\ \mbox{ for a.e. }x\in {\R}^k\right\}.
\]
 Precisely they satisfy
\begin{equation}\label{critical}
\left.\frac{d}{dt}{\mathcal{L}}^{1/2}(\Pi(u+t\varphi))\right|_{t=0}=0,\end{equation}
	where  $\varphi\in C^{1}_c(\R^k,\R^m)$, $\Pi:\mathcal{U}\to {\mathcal{N}}^n$ is any fixed $C^2$  projection, defined on some tubular neighborhood $\mathcal{U}$ of 
	${\mathcal{N}}^n$.
	
The homogeneous fractional Sobolev space  $ \dot H^{1/2}({\R}^k,{\R}^m)$ can be defined as follows
\[ \dot H^{1/2}(\R^k):=\bigg\{u\in L^2_{loc}(\R^k):\|u\|_{\dot H^{1/2}(\R^k)}^2:=\iint_{R^{2k}}\frac{|{u(x)-u(y)}|^2}{|{x-y}|^2}\,dx^k\,dy^k<\infty\bigg\}. \]
 The fractional Laplacian $(-\Delta)^{1/4}u$   can be  defined by
 
\[ (-\Delta)^{1/4}u:=\lim_{\epsilon\to 0}\mathcal{F}^{-1}[(\epsilon^2+4\pi^2|\xi|^2)^{1/4}\mathcal{F}u], \]
provided that the limit exists in $\mathcal{S}'(\R)$.\footnote{We denote respectively by ${\mathcal{S}}(\R)$   the spaces of (real or complex) Schwartz functions. Given a function $\varphi\in\mathcal{S}(\R)$, we denote either by
$\hat\varphi$ or by ${\mathcal{F}}\varphi$ the Fourier transform of $\varphi$, i.e.
\[ \hat\varphi(\xi)={\mathcal{F}}\varphi(\xi)=\int_{\R}v(x)e^{-2\pi i\xi x}\,dx. \]
}
\par 
We observe that  if  $u\in \dot H^{1/2}({\R}^k,{\R}^m)$,  then $(-\Delta)^{1/4} u$ exists, lies in $L^2(\R^k)$ and is given by
	\[ (-\Delta)^{1/4}u=\mathcal{F}^{-1}\left[\left(2\pi |{\xi}|\right)^{1/2}\widehat{u}\right], \]
(see for instance Lemma B.5 in \cite{DLP} and the references therein).

Weak  $1/2$-harmonic maps satisfy the  Euler-Lagrange equation  
\be
\label{int-I.3}
\nu(u)\wedge(-\Delta)^{1/2} u=0\quad \quad\mbox{in }{\mathcal D}'({\R}^k),
\ee
 where $\nu(z)$ is the Gauss Map  at $z\in {\mathcal{N}}^n$ taking values into the grassmannian $\ti{G}r_{m-n}({\R}^m)$ of oriented $m-n$ planes in ${\R}^m$ which
 is given by the oriented normal $m-n$-plane to $T_z {\mathcal{N}}^n\,.$\par
 One of the main result in \cite{DLR2} is the following Theorem
  \begin{Theorem}\label{holderhm}
Let  ${\mathcal{N}}^n$ be a   $C^2$ closed  submanifold of ${\R}^m$ and let $u\in \dot{H}^{1/2}(\R,{\mathcal{N}}^n)$ be a weak $1/2-$harmonic map into ${\mathcal{N}}^n$. Then  $u\in \bigcap_{0<\delta<1} C^{0,\delta}_{loc}(\R,{\mathcal{N}}^n).$ \hfill$\Box$
\end{Theorem}
Finally a   bootstrap argument leads to the following result  
(see \cite{DLP} for the details of this argument).
\begin{Theorem}\label{regul}
Let ${\mathcal{N}}^n\subset\R^m$ be a  $C^k$ closed submanifold of ${\R}^m$ , with $k\ge 2$, and let $u\in \dot{H}^{1/2}(\R,{\mathcal{N}}^n)$ be a weak $\frac{1}{2}$-harmonic. Then
\[ u\in
\bigcap_{0<\delta<1} C^{k-1,\delta}_{loc}(\R,{\mathcal{N}}^n). \]
In particular, if ${\mathcal{N}}^n$ is $C^\infty$ then $u\in C^\infty(\R,{\mathcal{N}}^n)$. \hfill$\Box$
\end{Theorem}
 We   remark that if   ${\mathcal{P}}_{-i}\colon S^1\setminus \{-i\}\to \R$, ${\mathcal{P}}_{-i}(\cos(\theta)+i \sin(\theta))=\frac{\cos(\theta)}{1+\sin(\theta)}$ is the classical stereographic projection whose inverse is given by
 \begin{equation}\label{formulaPi}
{\mathcal{P}}_{-i}^{-1}(x)=\frac{2x}{1+x^2}+i\left(-1+\frac{2}{1+x^2}\right),
\end{equation}
then the  following relation between the $1/2$-Laplacian in $\R$ and in $S^1$ holds:
 \begin{Proposition}  \label{proppull} Given $u:\R\to\R^m$, we set $v:=u\circ {\mathcal{P}}_{-i}:S^1\to\R^m$. Then $u\in L_{\frac{1}{2}}(\R)$\footnote{We recall that 
 $L_\frac{1}{2}(\R):=\left\{u\in L^1_{loc}(\R):\int_{\R}\frac{|u(x)|}{1+x^2}dx<\infty   \right\}$} if and only if $v\in L^1(S^1)$. In this case
\begin{equation}\label{eqlapv1}
\lapfr_{S^1} v(e^{i\theta})=\frac{((-\Delta)_{\R}^{\frac{1}{2}}u)({\mathcal{P}}_{-i}(e^{i\theta}))}{1+\sin\theta}  ~~\mbox{in $\mathcal{D}'(S^1\setminus\{-i\}).$}
\end{equation}
Observe that $(1+\sin(\theta))^{-1}=|{\mathcal{P}}_{-i}^{\prime}(\theta)|,$  hence we have
$$\int_{0}^{2\pi} \lapfr v(e^{i\theta})\,\varphi(e^{i\theta})\ d\theta=\int_{\R}\lapfr u(x)\ \varphi({\mathcal{P}}_{-i}^{-1}(x))\ dx\quad \text{for every }\varphi\in C^\infty_0(S^1\setminus\{-i\}).$$\end{Proposition}
For the proof of Proposition \ref{proppull}  we refer for instance to  \cite{DL1} or \cite{DLMR}.\par
From Proposition \ref{proppull}  it follows that $u\in {\dot  H}^{1/2}({\R})$ is a  $1/2$-harmonic map in $\R$  if and only if  $v:=u\circ {\mathcal{P}}_{-i}\in {\dot  H}^{1/2}({S^1})$  is a   $1/2$ harmonic map in $S^1$.\par
In the study of quantization properties of half-harmonic maps we established in \cite{DLLR}  new  {\bf Pohozaev-type  identities} for the half Laplacian in one dimension, that we are going to present below. \par
We first consider  the fundamental solution  $G$ of the fractional heat equation: 
 
 \begin{equation}\label{solfund1}
\left\{\begin{array}{cc}
\partial_t G+(-\Delta)^{1/2} G=0 & x\in\R,\,  t>0\\
G(0,x)=\delta_0 & t=0\,.
\end{array}\right.
\end{equation}
It is given by
  $$G(t,x)= \frac{1}{\pi} \frac{t}{x^2+t^2}.$$
The following equalities hold
$$\partial_t G=\frac{1}{\pi}\frac{x^2-t^2}{(t^2+x^2)^2},~~\partial_x G=-\frac{1}{\pi}\frac{2xt}{(t^2+x^2)^2}.$$

\begin{Theorem}\label{Pohozaev-R}{\bf [Pohozaev Identity in ${\mathbf{R}}$]}
Let $u\in W^{1,2}_{loc}(\R,\R^m)$ be  such that
\begin{equation}\label{Iharmequation1}
\frac{du}{dx}\cdot (-\Delta)^{1/2} u=0~~\mbox{a.e in $\R$.}
\end{equation}
Assume that for some $u_0\in R$
\begin{equation}\label{IcondPoh}
\int_{\R}|u-u_0 |dx <+\infty,~~\int_{\R}\left\vert \frac{du}{dx}(x)\right\vert\, dx <+\infty.
\end{equation}
Then the following identity holds
\begin{equation}\label{I-identityPR}
 \left|\int_{\R} \partial_tG(t,x)(u(x)-u_0)dx\right|^2=\left|\int_{\R} \partial_x G(t,x)  (u(x)-u_0)dx\right|^2 \hbox{ for all } t\in \R.~\Box
 \end{equation}
\end{Theorem}

\par
  \bigskip
We get an analogous formula in $S^1$.
By  identifying  $S^1$ with $[-\pi,\pi)$ we consider the following problem
\begin{equation}\label{solfund3}
\left\{\begin{array}{cc}
\partial_t F+(-\Delta)^{1/2} F=0 & ~~\theta\in[-\pi,\pi),\; t>0\\
F(0,\theta)=\delta_0(x) & \theta\in [-\pi,\pi].
\end{array}\right.
\end{equation}
  The solution of \rec{solfund3} is given by
$$
F(\theta,t)=\frac{1}{2\pi}\sum_{n=-\infty}^{+\infty} e^{-t|n|}e^{in\theta}=\frac{e^{2t}-1}{e^{2t}-2e^{t}\cos(\theta)+1}.$$\par
In this case we have
$${\partial_t F} (t,\theta)=-2e^t\frac{e^{2t}\cos(\theta)-2e^{t}+\cos(\theta)}{(e^{2t}-2e^{t}\cos(\theta)+1)^2}$$
and
$$
{\partial_{\theta} F} (t,\theta) =-2e^{t}\frac{ \sin(\theta)(e^{2t}-1)}{(e^{2t}-2e^{t}\cos(\theta)+1)^2}.$$
Then the following holds
\begin{Theorem}\label{Pohozaev-S^1}
{\bf [Pohozaev Identity on ${\mathbf S^1}$]} 
Let $u\in   W^{1,2}(S^1,\R^m)$ be such that
\begin{equation}\label{harmequation2}
\frac{\partial u}{\partial \theta}\cdot (-\Delta)^{1/2} u=0~~\mbox{a.e in $S^1$.}\end{equation}
Then the following identity holds
\begin{equation}\label{pohos1}
 \left|\int_{S^1}u(z)\partial_t F(z) \, d\theta\right|^2= \left|\int_{S^1}u(z)\partial_{\theta}F(z) \, d\theta\right|^2.
\end{equation}
From \eqref{pohos1} one deduces in particular that 
\begin{equation}\label{poho-s1}
 \left|\int_{0}^{2\pi}u(e^{i\theta}) \cos(\theta)  \, d\theta\right|^2= \left|\int_{0}^{2\pi}u(e^{i\theta}) \sin(\theta)  \, d\theta\right|^2.~~~\Box
\end{equation}
\end{Theorem}\medskip
  For the proof of Theorem \ref{Pohozaev-R} and Theorem \ref{Pohozaev-S^1} we refer the reader to \cite{DLLR}.\par
         \par
   We could have solved \eqref{solfund1} by requiring $G(0,x)=\delta_{x_0}$, with $x_0\in\R$ and we would have obtained as many corresponding Pohozaev-type formulas.\par
We observe that  if   $u $ is a  smooth critical point of 
  \eqref{fracenergy} in $\R$ then it is  {\em stationary } as well, namely it is critical with respect to the variation of the domain:
  \begin{equation}\label{statio}
\left(  \frac{d}{da} \int_{\R}| (-\Delta)^{1/4}\Pi(u(x+aX(x))|^2 dx\right)_{\big |_{a=0}}=\left(  \frac{d}{da} \int_{\R}| (-\Delta)^{1/4}(u(x+aX(x))|^2 dx\right)_{\big |_{a=0}}=0
 \end{equation}
 where $X\colon \R^2\to \R^2$ is a $C^1_c(R^2)$ vector field.   \par
  Actually any variation the form $\displaystyle{u(x+aX(x))=u(x)+a\frac{du(x)}{dx} X(x)+o(a^2)}$ can be interpreted  
as being a variation in the target with $\displaystyle{\varphi(x)= \frac{du(x)}{dx} X(x)}.$

From   \eqref{statio} we get the so-called equation of stationarity:
  \begin{equation}\label{eqstaR}
0=\int_{\R}[(-\Delta)^{1/2}(u(x+aX(x))\cdot \frac{d}{da}(u(x+aX(x)))]_{\big  |_{a=0}}dx =\int_{\R}(-\Delta)^{1/2}(u(x))\cdot \frac{du(x)}{dx} X(x)dx. \nonumber
\end{equation}
By the arbitrariness of $ X$  and the smoothness of $u$ from \eqref{eqstaR} we deduce that

\begin{equation}\label{eqstat1}
(-\Delta)^{1/2}u (x)\cdot \frac{du }{dx}(x)=0~~\mbox{ $x\in \R.$}
\end{equation}
In an analogous way    if $u$ is a smooth  critical point of     the fractional energy  \eqref{fracenergy}  in $S^1$,  it also satisfies
 \begin{equation}\label{variatstat}
\left( \frac{d}{da}\int_{S^1}| (-\Delta)^{1/4}(u(z+aX(z)))|^2d\sigma(z)\right) _{\big |_{a=0}}=0
 \end{equation}
  where $X\colon S^1\to \R^2$ is $C^1(S^1)$ vector field. 
    \begin{equation}\label{eqstat2}  (-\Delta)^{1/2}(u(z))\cdot\partial_{z}u(z) =0~~\mbox{   $z\in S^1$}.\end{equation}

   Therefore the assumptions of Theorem \ref{Pohozaev-R} and Theorem \ref{Pohozaev-S^1} are satisfied by sufficiently smooth  $1/2$-harmonic maps.\par

We have now to give some explanations why these identities belong to the {\it Pohozaev identities} family.  These identities are produced by the conformal invariance of the highest
order derivative term in the Lagrangian from which the Euler Lagrange is issued.  For instance the Dirichlet energy
\begin{equation}\label{DE}
{\mathcal{L}}(u)=\int_{\R^2} |\nabla u|^2 dx^2
\end{equation}
is conformal invariant in $2$-D, whereas the following fractional energy
\begin{equation}\label{FE}
{\mathcal{L}}_{\R}^{1/2}(u)=\int_{\R} |(-\Delta)^{1/4}u|^2\ dx
\end{equation}
is conformal invariant  in $1$-D. The infinitesimal perturbations issued from the dilations produce in \rec{DE} and \rec{FE} respectively   the following infinitesimal variations of these highest order terms 
\[
  \sum_{i=1}^2 x_i\,\frac{\p u}{\p x_i}\cdot\Delta u\quad\mbox{ in $2$-D}\quad\mbox{ and }\quad x\frac{d u}{dx}\cdot(-\Delta)^{1/2}u\quad\mbox{ in  $1$-D}
\]
Such kind of  perturbations  play an important role in establishing Pohozaev-type identities.\par

 In two dimensions, integrating the identity \eqref{eqstat1} on a ball $B(x_0,r)$  $(x_0\in R^2,r>0)$  gives the following {\em  balancing law} between the {\em radial} part and the {\em angular} part of the energy
classically known as {\bf Pohozaev identity}.
\begin{Theorem}\label{th-I.6} 
Let $u\in W^{2,2}_{loc} ({\R}^2,\R^m)$ such that
\begin{equation}\label{IharmequationR2}
 \frac{\partial u}{\partial x_i }\cdot \Delta u=0~~\mbox{a.e in $B(0,1)$}
\end{equation}
for $i=1,2$.
 Then it holds
\begin{equation}
\label{PI}
\int_{\partial B(x_0,r)} \left\vert \frac{1}{r}\frac{\partial u}{\partial\theta}\right\vert^2 d\theta=\int_{\partial B(x_0,r)}\left\vert\frac{\partial u}{\partial r}\right\vert^2 d\theta
\ee
for all $r\in [0,1].$
\end{Theorem}
 
\begin{figure}
 
 \begin{center}
\includegraphics[width=8cm]{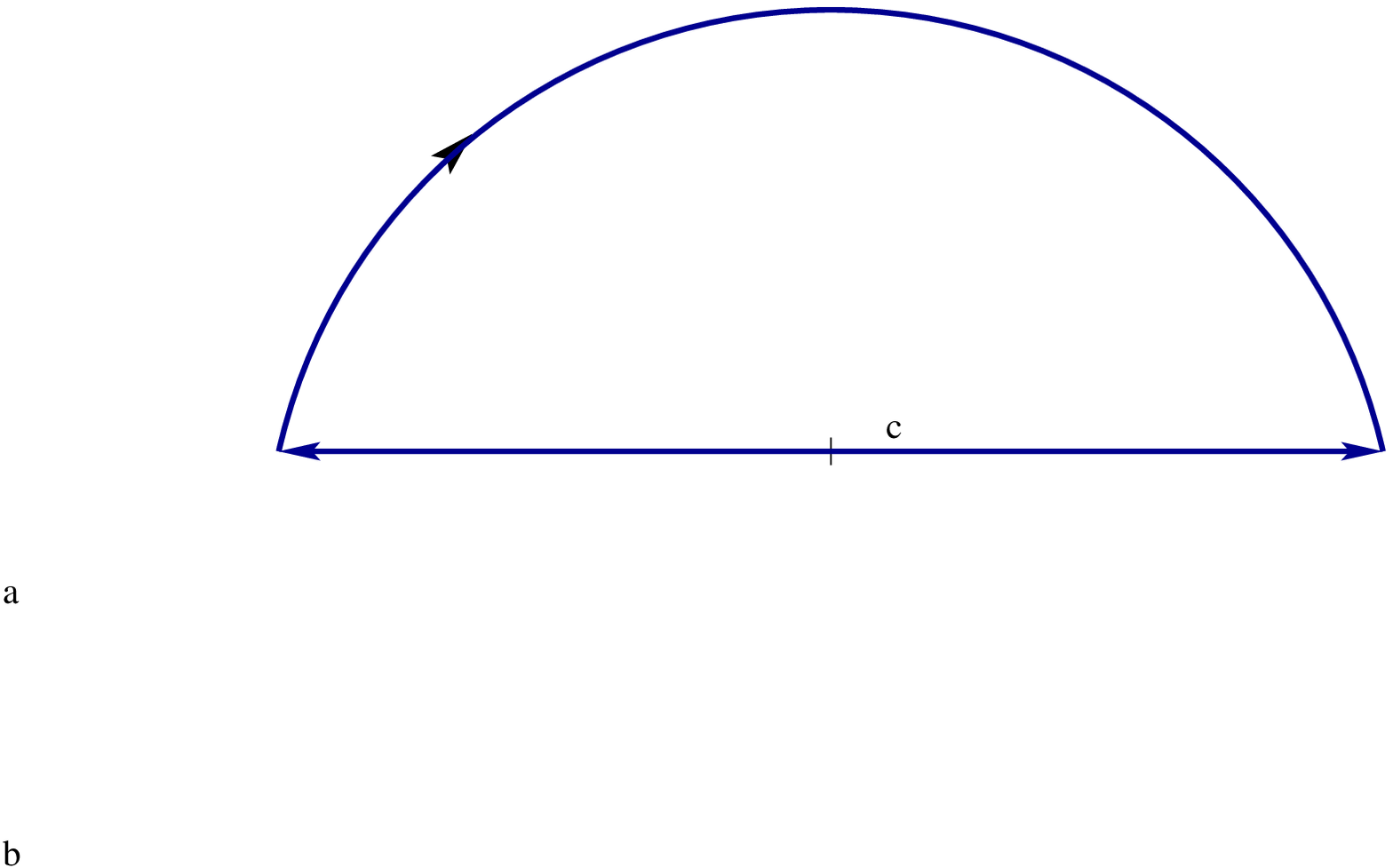}
\end{center}
 \vskip1cm
  \begin{center}
\caption{Link between the symmetric and antisymmetric part of $u$ and the integral of the radial and tangential derivative of any extension $\tilde u$ of $u$ on upper half plane $R^2_+$}\label{fig1}
\end{center}
\end{figure}

In $1$ dimension one might wonder what corresponds to the 2 dimensional dichotomy between {\em radial} and {\em
 angular} parts. Figure \ref{fig1}  is intended to illustrate the following correspondence
of dichotomies respectively in $1$ and $2$ dimensions.
 
\psfrag{a}{$  \displaystyle{\int_{\partial B(0,r)\cap R^2_+}\frac{\partial \tilde u}{\partial \theta}d\theta=(u(r)-u(0))-(u(-r)-u(0))=(u-u(0))^{-}(r)}$}
\psfrag{b}{$ \displaystyle{\int_{\partial B(0,r)\cap R^2_+}\frac{\partial \tilde u}{\partial r}dr =u(r)+u(-r)-2u(0)=(u-u(0))^{+}(r)}$}
\psfrag{c}{$0$}
 
\[
\begin{array}{ccc}
\mbox{$2$-D}& \longleftrightarrow & \mbox{ $1$-D}\\[3mm]
\ds\mbox{ {\bf radial} : } \quad\frac{\p u}{\p r}& \longleftrightarrow & \mbox{ {\bf symmetric} part of }u\quad:\quad u^+(x):= \frac{u(x)+u(-x)}{2}\\[3mm]
\ds\mbox{{\bf  angular} : }\quad \frac{\p u}{\p \theta}& \longleftrightarrow & \mbox{ {\bf antisymmetric} part of }u\quad:\quad  u^-(x):= \frac{u(x)-u(-x)}{2}
\end{array}
\]
In this note we make the observation  that by exploiting the  invariance of the equation \eqref{harmequation2} with respect to the trace of M\"obius transformations of the disk in $\R^2$  of the form $M_ {\alpha,a}(z):= e^{i\alpha}\frac{z-{a}}{1- az}$, $\alpha\in\R, a\in (-1,1)$ \footnote {We recall that
since  $M_ {\alpha,a}(z)$ is conformal  with $M'_ {\alpha,a}(z)\ne 0$ we have
\begin{equation}
(-\Delta)^{1/2}(u\circ M_ {\alpha,a}(z))=(-\Delta)^{1/2}u_a= e^{\lambda_{\alpha,a}} ((-\Delta)^{1/2}u)\circ M_ {\alpha,a}(z),
\end{equation}
where $\lambda_{\alpha,a}(z)=\log(|\frac{\partial M_ {\alpha,a}}{\partial\theta}(z)|), z\in S^1$}
we can derive 
from \eqref{poho-s1} a countable family of  relations  involving the Fourier coefficients of solutions of \eqref{harmequation2}. This fact has been already announced in the paper \cite{DLLR}.
We heard  that the proof of this property has been recently obtained also in the work of preparation \cite{CMM} by using a different approach.\par
 Given $u\colon S^1\to\R^m$ we define its  Fourier coefficients
 for every $k\ge 0$ :
 \[
\left\{\begin{array}{l}
\ds a_k:=\frac{1}{2\pi}\int_0^{2\pi}u(e^{i\theta})\, \cos k\theta\ d\theta\\[3mm]
\ds b_k=\frac{1}{2\pi}\int_0^{2\pi}u(e^{i\theta})\, \sin k\theta\ d\theta.
\end{array}
\right.
\]

The following result holds.
\begin{Proposition} 
 \label{Pohozaev2-S^1}
{\bf [Relations of  the Fourier coefficients on ${\mathbf S^1}$]}
Let   $u\in W^{1,2}(S^1,\R^m)$ satisfy 
\eqref{harmequation2}.
Then for every $n\ge 2$ it holds
\begin{equation}\label{identiFC0}
\sum_{k=1}^{n-1}(n-k)k (a_ka_{n-k}-b_kb_{n-k})=0 \end{equation}
and
 \begin{equation}\label{identiFC02}
 \sum_{k=1}^{n-1}(n-k)k(a_kb_{n-k}+b_ka_{n-k}) =0.~~~~\Box
\end{equation}
  \end{Proposition}
 
We conclude this   introduction by mentioning  that in the   paper \cite{ROS} the authors obtains a different Pohozaev identity for bounded weak solutions to the following problem
\begin{equation}\label{nlpl}
\left\{\begin{array}{cc}
(-\Delta)^{s}u=f(u)&~\mbox{in $\Omega$}
\\
u=0&~\mbox{in $\R^n\setminus\Omega$}
\end{array}\right.
\end{equation}
where $s\in(0,1)$ and $\Omega\subset\R^n$ is a bounded domain.
As a consequence of their Pohozaev identity they get nonexistence results for problem \eqref{nlpl} with supercritical nonlinearitis in star-shaped domains.\par
This paper is organized as follows. In section \ref{sect2} we prove Proposition \ref{Pohozaev2-S^1} and in section \ref{sect3} we provide      infinite many   Pohozaev formulas for stationary harmonic maps  in $2$-D in  correspondence to conformal  vector fields in $\C$  generated by the holomorphic functions. In particular we will prove a generalization of Theorem \ref{th-I.6}.  This can be  useful for different purposes.

\section{Proof of Proposition \ref{Pohozaev2-S^1}. }\label{sect2}
 From Theorem \ref{Pohozaev-S^1} it follows that $u$ satisfies   in particular  
\begin{equation}\label{poho-s1bis}
 \left|\int_{0}^{2\pi}u(e^{i\theta}) \cos(\theta)  \, d\theta\right|^2= \left|\int_{0}^{2\pi}u(e^{i\theta}) \sin(\theta) \, d\theta\right|^2
 \end{equation}

 \par

We can rewrite  (\ref{poho-s1bis}) as follows
  \be
\label{poho-s3}
\left|\int_0^{2\pi}u(e^{i\theta})\, \Re(d e^{i\theta})\right|^2=\left|\int_0^{2\pi}u(e^{i\theta})\, \Im(d e^{i\theta})\right|^2.
\ee
  
Given $a\in \R$ with $|a|<1$ and $\alpha\in\R$  we consider the  M\"obius map  $M_ {\alpha,a}(z):= e^{i\alpha}\frac{z-{a}}{1- az}$ and we define  $$u_{a,\alpha}(e^{i\theta}):=u\circ M_ {\alpha,a}(z).$$
Since  the condition \eqref{harmequation2} is invariant with respect to M\"obius  transformations for every $\alpha\in\R$ and for every $a\in(-1,1)$ we get
\be
\label{poho-s4}
\left|\int_0^{2\pi}u\left(e^{i\alpha}\frac{z-{a}}{1- az}\right)\, \Re(d e^{i\theta})\right|^2=\left|\int_0^{2\pi}u\left(e^{i\alpha}\frac{z-{a}}{1- az}\right)\, \Im(d e^{i\theta})\right|^2. 
\ee
or equivalently 
 \be
\label{poho-s5}
\left|\Re\left(\int_0^{2\pi}u\left(e^{i\alpha}\frac{e^{i\theta}-{a}}{1- ae^{i\theta}}\right)\, d e^{i\theta}\right)\right|^2=\left|\Im\left(\int_0^{2\pi}u\left(e^{i\alpha}\frac{e^{i\theta}-{a}}{1- ae^{i\theta}}\right)\, d e^{i\theta}\right)\right|^2. 
\ee
We set $$
e^{i\varphi}:=e^{i\alpha}\frac{e^{i\theta}-a}{1- ae^{i\theta}},$$
which implies that
\be\label{subst1}
e^{i\theta}=\frac{e^{i(\varphi-\alpha)}+a}{1+ae^{i(\varphi-\alpha)}}
\ee
\be\label{subst2}
d(e^{i\theta})=\frac{1-a^2}{(1+ae^{i(\varphi-\alpha)})^2}d(e^{i(\varphi-\alpha)})\ee

By plugging (\ref{subst1}) and (\ref{subst2}) into (\ref{poho-s5}) and dividing by  $(1-a^2)$  we get
 \be
\label{poho-s6}
\left|\Re\left(\int_0^{2\pi}u(e^{i\varphi})\, \frac{e^{-i\alpha} }{(1+ae^{i(\varphi-\alpha)})^2}d(e^{i\varphi})\right)\right|^2=\left|\Im\left(\int_0^{2\pi}u(e^{i\varphi})\, \frac{e^{-i\alpha} }{(1+ae^{i(\varphi-\alpha)})^2}d(e^{i\varphi}\right)\right|^2. 
\ee
Observe that for all $|z|<1$ we have
$$
\frac{z}{(1+z)^2}=\sum_{n=1}^{\infty}n(-1)^{n-1}z^n$$
In particular
\be
 \frac{e^{i(\varphi-\alpha)}}{(1+ae^{i(\varphi-\alpha)})^2}= \sum_{n=1}^{\infty}n(-1)^{n-1}a^{n-1}e^{in(\varphi-\alpha)}
 \ee
 and
 \be
 \Re\left(\frac{e^{i(\varphi-\alpha)}}{(1+ae^{i(\varphi-\alpha)})^2}\right)= \sum_{n=1}^{\infty}n(-1)^{n-1}a^{n-1}\cos(n(\varphi-\alpha))
 \ee
  \be
 \Im\left(\frac{e^{i(\varphi-\alpha)}}{(1+ae^{i(\varphi-\alpha)})^2}\right)= \sum_{n=1}^{\infty}n(-1)^{n-1}a^{n-1}\sin(n(\varphi-\alpha))
 \ee
 We can write
  \begin{align}
\label{poho-s7}
&\left|\Re\left(\int_0^{2\pi}u(e^{i\varphi})\, \frac{e^{i(\varphi-\alpha)}}{(1+ae^{i(\varphi-\alpha)})^2}d\varphi\right)\right|^2\\
&=\sum_{n=1}^{\infty}(-1)^{n-1}a^{n-1}\sum_{k=1}^{n-1}(n-k)k\left(\int_0^{2\pi}u(e^{i\varphi})\cos(k(\varphi-\alpha))d\varphi\right)
\left(\int_0^{2\pi}u(e^{i\varphi})\cos((n-k)(\varphi-\alpha))d\varphi\right)\nonumber
\end{align}
and
\begin{align}
\label{poho-s8}
 &\left|\Im\left(\int_0^{2\pi}u(e^{i\varphi})\, \frac{e^{i(\varphi-\alpha)}}{(1+ae^{i(\varphi-\alpha)})^2}d\varphi\right)\right|^2\\
&=\sum_{n=1}^{\infty}(-1)^na^{n-1}\sum_{k=1}^{n-1}(n-k)k\left(\int_0^{2\pi}u(e^{i\varphi})\sin(k(\varphi-\alpha))d\varphi\right)
\left(\int_0^{2\pi}u(e^{i\varphi})\sin((n-k)(\varphi-\alpha))d\varphi\right)\nonumber
\end{align}
The identity (\ref{poho-s6}) and the relations (\ref{poho-s7}), (\ref{poho-s8}) imply that for every $n\ge 2$ we obtain the following identities
 \begin{align}\label{identiFC}
&\sum_{k=1}^{n-1}(n-k)k\left(\int_0^{2\pi}u(e^{i\varphi})\cos(k(\varphi-\alpha))d\varphi\right)\left(\int_0^{2\pi}u(e^{i\varphi})\cos((n-k)(\varphi-\alpha))d\varphi\right)\\
&=\sum_{k=1}^{n-1}(n-k)k\left(\int_0^{2\pi}u(e^{i\varphi})\sin(k(\varphi-\alpha))d\varphi\right)\left(\int_0^{2\pi}u(e^{i\varphi})\sin((n-k)(\varphi-\alpha))d\varphi\right).\nonumber
\end{align}
From \eqref{identiFC} we can deduce a countable family of relations between the Fourier coefficients of the map $u.$
Precisely if we set  for every $n\ge 1$
\[
\left\{\begin{array}{l}
\ds a_n:=\frac{1}{2\pi}\int_0^{2\pi}u(e^{\theta})\, \cos n\theta\ d\theta\\[3mm]
\ds b_n=\frac{1}{2\pi}\int_0^{2\pi}u(e^{\theta})\, \sin n\theta\ d\theta,
\end{array}
\right.
\]
  we get 
 \begin{align}\label{identiFC2}
&\sum_{k=1}^{n-1}(n-k)k\left[\left(\cos(k\alpha)a_k+\sin (k\alpha)b_k\right)\left(\cos((n-k)\alpha)a_{n-k}+\sin ((n-k)\alpha)b_{n-k}\right)\right.\nonumber\\
&~~~~~-\left.\left(\cos(k\alpha)b_k-\sin (k\alpha)a_k\right)\left(\cos((n-k)\alpha)b_{n-k}-\sin ((n-k)\alpha)a_{n-k}\right)\right]=0
 \end{align}
 The identity \eqref{identiFC2} can be rewritten as follows
 \begin{align}\label{identiFC3}
&\cos(n\alpha)(\sum_{k=1}^{n-1}(n-k)k (a_ka_{n-k}-b_kb_{n-k})) +\sin(n\alpha)(\sum_{k=1}^{n-1}(n-k)k(a_kb_{n-k}+b_ka_{n-k}) )=0.
\end{align}
The relation \eqref{identiFC3} yields \eqref{identiFC0} and \eqref{identiFC02} because of the linear dependence of $\cos(n\alpha)$ and $\sin(n\alpha)$.\par

We observe that
  for $n=2$ we obtain:
\begin{equation}\label{firstFC}
(|a_1|^2-|b_1|^2)\cos(2\alpha)-2a_1\cdot b_1 \sin(2\alpha)=0.
 \end{equation}
 
 Since $\alpha\in\R$ is arbitrary we get
 
\[
\left\{\begin{array}{l}
\ds |a_1|=|b_1|\\[3mm]
\ds a_1\cdot b_1=0
\end{array}
\right.
\]

If $n=3$ we get
\begin{equation}\label{firstFC2}
4( a_1\cdot a_2-b_1\cdot b_2)\cos(3\alpha)-4(a_1\cdot b_2+b_1\cdot a_2) \sin(3\alpha)=0.
 \end{equation}
 The relation \eqref{firstFC2} gives
\[
\left\{\begin{array}{l}
\ds a_1\cdot a_2=b_1\cdot b_2\\[3mm]
\ds a_1\cdot b_2=-a_2\cdot b_1.
\end{array}
\right.
\]
 
If $n=4$ we get
 
\[
\left\{\begin{array}{l}
\ds |a_2|^2-|b_2|^2=\frac{3}{2}( b_1\cdot b_3-a_1\cdot a_3)\\[3mm]
\ds a_2\cdot b_2=-\frac{3}{4}(a_1\cdot b_3+b_1\cdot a_3).
\end{array}
\right.
\]

We can conclude the proof.~~~~~$\Box$

 \section{Pohozaev Identities for the Laplacian in $\R^2$}\label{sect3}
    In this section we derive  Pohozaev identities in $2$-D in the same spirit of the previous section.
    In the following Theorem \ref{Poho0} we combine ideas from \cite{PR} and \cite{Stru}.
    Precisely we multiply the equation satisfied for instance by sufficiently smooth harmonic maps      by   the fundamental solution of the heat equation and a holomorphic vector field
    $X\colon\C\to\C.$\par
    We mention that    the use  of the fundamental solution to get Pohozaev-type identities and monotonicity formulas 
    has been performed in \cite{Stru} to study the heat flow.  In Chapter 9 of \cite{PR} the authors derived in the context of Ginzurg-Landau equation generalized Pohozaev identities for the so-called {\em $\rho$- conformal   vector fields}
    $X=(X^1,\ldots,X^n)$, where $\rho$ is a given function defined in a $2$ dimensional domain. In the case $\rho\equiv 1$    then the   $\rho$- conformal   vector fields are exactly 
    conformal  vector fields   generated by the holomorphic functions.\par
We recall that the fundamental solution of the heat equation
\begin{equation}\label{solfundlap}
\left\{\begin{array}{cc}
\partial_t G+(-\Delta) G=0 & t>0\\
G(0,x)=\delta_{x_0} & t=0\,.
\end{array}\right.
\end{equation}
is given by $G(x,t)=(4\pi t)^{-1/2}e^{-\frac{|x-x_0|^2}{4t}}.$
   \begin{Theorem}{[Pohozev in  $\R^2$]}\label{Poho0}
Let $u\in W^{2,2}_{loc}(\R^2,\R^m)$ such that
\begin{equation}\label{harmequation0}
 \frac{\partial u}{\partial x_i}\cdot \Delta u =0~~\mbox{a.e in $\R^2$}
\end{equation}
for $i=1,2.$
Assume that
\begin{equation}\label{cond1}
 \int_{\R^2}|\nabla u(x)|^2dx<+\infty.
\end{equation}
Then for all $x_0\in\R^2$,  $t>0$ and every $X=X_1+i X_2\colon \C\to \C$ holomorphic function the  following identity holds
\begin{equation}
  2  \iint_{R^2}e^{-{\frac{|x-x_0|^2}{4t}}}|x-x_0|  \left(\frac{\partial u}{\partial \nu}\cdot  \frac{\partial u}{\partial X}  \right)dx=   \iint_{R^2}e^{-{\frac{|x-x_0|^2}{4t}}}\left(( x-x_0)\cdot  X\right) |\nabla u|^2 dx.
  \end{equation}
  In the particular case where $X=x-x_0$  with  $x_0\in\R^2$   
  then for all $t>0$   the following identity holds
\begin{equation}\label{idpohozaevR2}
  \iint_{R^2}e^{-{\frac{|x-x_0|^2}{4t}}}|x-x_0|^2\left|\frac{\partial u}{\partial \nu} \right|^2dx^2=  \iint_{R^2}e^{-{\frac{|x-x_0|^2}{4t}}}\left|\frac{\partial u}{\partial \theta}\right|^2 dx^2.
  \end{equation}
\end{Theorem}
 {\bf Proof.}  We multiply the equation \rec{harmequation0} by $X_i e^{-{\frac{|x-x_0|^2}{4t}}}$ and we integrate\footnote{We use the Einstein summation convention}
\begin{eqnarray}\label{estpoh0}
0&=&\sum_{k,i=1}^2\iint_{\R^2}  e^{-{\frac{|x-x_0|^2}{4t}}}X_i\frac {\partial u}{\partial x_i}\frac{\partial^2 u}{\partial x_k^2}dx\nonumber\\
&=&\sum_{k,i=1}^2\frac{1}{2t}\iint_{\R^2}  e^{-{\frac{|x-x_0|^2}{4t}}}(x_k-x_{k,0})X_i\frac {\partial u}{\partial x_i}\frac{\partial u}{\partial x_k}dx\nonumber\\
&-&\sum_{k,i=1}^2\iint_{\R^2}  e^{-{\frac{|x-x_0|^2}{4t}}} \frac{\partial X_i}{\partial x_k}\frac {\partial u}{\partial x_i}\frac{\partial u}{\partial x_k}dx\nonumber\\
&-&\sum_{k,i=1}^2\iint_{\R^2}  e^{-{\frac{|x-x_0|^2}{4t}}} X_i \frac{\partial^2 u}{\partial x_k\partial x_i}\frac{\partial u}{\partial x_k}dx\nonumber\\
&=&\frac{1}{2t}\iint_{\R^2}e^{-{\frac{|x-x_0|^2}{4t}}} (X\cdot\nabla u)\frac{\partial u}{\partial \nu}|x-x_0|dx\\
&-&\iint_{\R^2}e^{-{\frac{|x-x_0|^2}{4t}}}\left[\underbrace{\frac{\partial X_1}{\partial x_1}}_{\frac{\partial X_1}{\partial x_1}=\frac{\partial X_2}{\partial x_2}}|\nabla u|^2+(\underbrace{\frac{\partial X_1}{\partial x_2}+\frac{\partial X_2}{\partial x_1}}_{=0})(\frac{\partial u}{\partial x_1}\frac{\partial u}{\partial x_2})\right]\nonumber\\
&-&\iint_{\R^2}e^{-{\frac{|x-x_0|^2}{4t}}}X_i\frac{1}{2}\partial_{x_i}|\nabla u|^2 dx\nonumber\\
&=&\frac{1}{2t}\iint_{\R^2}e^{-{\frac{|x-x_0|^2}{4t}}} (X\cdot\nabla u)\frac{\partial u}{\partial \nu}|x-x_0|dx\nonumber\\
&-&\frac{1}{4t}\iint_{\R^2}e^{-{\frac{|x-x_0|^2}{4t}}}(x-x_0)\cdot X |\nabla u|^2 dx\nonumber\\
&-&\iint_{\R^2}e^{-{\frac{|x-x_0|^2}{4t}}} \frac{\partial X_1}{\partial x_1} |\nabla u|^2\nonumber\\
&+&\frac{1}{2}\sum_{i=1}^2\iint_{\R^2} e^{-{\frac{|x-x_0|^2}{4t}}}  \frac{\partial X_i}{\partial x_i} |\nabla u|^2 dx \nonumber
\end{eqnarray}
By using the fact that 
$$\frac{1}{2}\sum_{i=1}^2\iint_{\R^2} e^{-{\frac{|x-x_0|^2}{4t}}}  \frac{\partial X_i}{\partial x_i} |\nabla u|^2 dx=\iint_{\R^2}e^{-{\frac{|x-x_0|^2}{4t}}} \frac{\partial X_1}{\partial x_1} |\nabla u|^2$$
we obtain from \eqref{estpoh0} that
\begin{equation}\label{estpoh1}
2\iint_{\R^2}e^{-{\frac{|x-x_0|^2}{4t}}} (X\cdot\nabla u)\frac{\partial u}{\partial \nu}|x-x_0|dx=
 \iint_{\R^2}e^{-{\frac{|x-x_0|^2}{4t}}}(x-x_0)\cdot X |\nabla u|^2 dx.
\end{equation}
 In particular if $X=(x-x_0)$ by using that $\nabla u=(\frac{\partial u}{\partial\nu},|x-x_0|^{-1}\frac{\partial u}{\partial\theta})$, from \rec{estpoh1} we get the identity
  
\begin{equation} 
  \iint_{R^2}e^{-{\frac{|x-x_0|^2}{4t}}}|x-x_0|^2\left|\frac{\partial u}{\partial \nu} \right|^2dx=  \iint_{R^2}e^{-{\frac{|x-x_0|^2}{4t}}}\left|\frac{\partial u}{\partial \theta}\right|^2 dx
  \end{equation}
  and we conclude.~\hfill$\Box$\par
  \medskip
  In   Theorem \ref{Pohobis} we get infinite many    Pohozaev identities over balls in correspondence to  holomorphic vector fields   $X=X_1+i X_2\colon \C\to \C$ for maps   $u\in W^{2,2}_{loc}(\R^2,\R^m)$ satisfying \eqref{harmequation0}

  \begin{Theorem}{[Pohozev in  $\R^2$- Ball Case]}\label{Pohobis}
Let $u\in W^{2,2}_{loc}(\R^2,\R^m)$ such that
\begin{equation}\label{harmequationball}
 \frac{\partial u}{\partial x_i}\cdot \Delta u=0 ~~\mbox{a.e in $\R^2$}
\end{equation}
for $i=1,2.$
 Then for all $x_0\in\R^2$,  $r>0$ and every $X=X_1+i X_2\colon \C\to \C$ holomorphic function the  following identity holds
\begin{equation}
  2  \int_{\partial B(x_0,r)}    \frac{\partial u}{\partial \nu} \,\nabla u\cdot X dx=   \int_{\partial B(x_0,r)}  X\cdot\nu |\nabla u|^2 dx
  \end{equation}
  In the particular case where $X=x-x_0$  
 with  $x_0\in\R^2$,   then for all $r>0$  the  following identity holds
 \begin{equation}
  2  \int_{\partial B(x_0,r)}    \left|\frac{\partial u}{\partial \nu}\right|^2 d\sigma =   \int_{\partial B(x_0,r)} \left|\nabla u\right|^2 d\sigma.
  \end{equation}
or
 \begin{equation}
   \int_{\partial B(x_0,r)}    \left| \frac{\partial u}{\partial \nu}\right|^2 d\sigma = \frac{1}{r^2}  \int_{\partial B(x_0,r)}\left|\frac{\partial u}{\partial \theta}\right|^2 d\sigma 
  \end{equation}
 \end{Theorem}
 {\bf Proof.}  We multiply the equation \rec{harmequationball} by $X_i $ and we integrate over $B(x_0,r)$:
 \begin{eqnarray}\label{estphoballs}
 0&=&
 \int_{B(x_0,r)} X_i\frac {\partial u}{\partial x_i}\frac{\partial^2 u}{\partial x_k^2}dx\nonumber\\
&=& \int_{B(x_0,r)} \frac{\partial}{\partial x_k}\left(X_i \frac {\partial u}{\partial x_i}\frac {\partial u}{\partial x_k}\right)dx\\
&-& \int_{B(x_0,r)} \frac{\partial X_i}{\partial x_k} \frac {\partial u}{\partial x_i}\frac {\partial u}{\partial x_k}dx\nonumber\\
&-& \frac{1}{2}\int_{B(x_0,r)} X_i \frac {\partial}{\partial x_i}|\nabla u|^2 dx\nonumber\\
&=&\int_{\partial B(x_0,r)} (X \cdot \nabla u )(\frac{\partial u}{\partial \nu}) d\sigma\nonumber\\
&-&\frac{1}{2r}\int_{\partial B(x_0,r)} X\cdot (x-x_0)  |\nabla u|^2 dx\nonumber\\
&+&\frac{1}{2}\int_{B(x_0,r)} \frac {\partial X_i}{\partial x_i}  |\nabla u|^2 dx- \int_{B(x_0,r)} \frac{\partial X_i}{\partial x_k} \frac {\partial u}{\partial x_i}\frac {\partial u}{\partial x_k}dx.\nonumber
\end{eqnarray}
By using the Cauchy Riemann equations one deduces that the last term in \eqref{estphoballs} is zero and therefore

$$\int_{\partial B(x_0,r)} (X \cdot \nabla u )( \nabla u\cdot (x-x_0) )d\sigma=\frac{1}{2}\int_{\partial B(x_0,r)} X\cdot (x-x_0)  |\nabla u|^2 d\sigma.
$$
and we conclude.\hfill$\Box$

\end{document}